\documentclass[12pt]{article}
\usepackage{latexsym}
\usepackage{amsfonts}

\usepackage[active]{srcltx} 

\newcommand{\ba}{\begin{array}}\newcommand{\ea}{\end{array}}
\newcommand{\sym}{S}
\newcommand{\ns}{\rm}

\renewcommand{\Bbb}{\mathbb}
\newcommand{\nse}{\kern-3pt\ns$=$}\newcommand{\qd}{\hfill$\Box$\medbreak}

\newcommand{\Z}{\mbox{$\Bbb Z$}}
\newcommand{\C}{\mbox{$\Bbb C$}}
\newcommand{\noi}{\noindent}
\newcommand{\ext}{\raise1pt\hbox{$\ts\bigwedge$}}
\renewcommand{\sym}{S}
\newcommand{\ts}{\textstyle}
\newcommand{\rf}[1]{(\ref{#1})}\renewcommand{\c}[1]{\cite{#1}}
\newcommand{\chii}{\raise2pt\hbox{$\chi$}}

\newcommand{\la}{\lambda}

\newcommand{\Fg}{\mbox{${\cal F}\kern-2pt_g$}}

\newcommand{\Mg}{\mbox{${\cal M}\kern-2pt_g$}}
\newcommand{\Ng}{\mbox{${\cal N}\kern-2pt_g$}}
\newcommand{\V}{V\kern-1pt}

\newcommand{\Gg}{\mbox{${\cal G}\kern-2pt_g$}}
\newcommand{\cir}{\raise1.6pt\hbox{\footnotesize$\circ$}}

\newcommand{\ind}{\mbox{\ns ind}}

\newcommand{\mod}{\mbox{\ns mod\,\,\,}}
\newcommand{\sign}{\mbox{${\rm sign}$}}
\newcommand{\Res}[2]{\hbox{\ns

 Res}\kern-16pt\lower5pt\hbox{\footnotesize$_{#1}$}\kern2pt\left[#2\right]}
\newcommand{\qk}{quaternion-K\"ahler\kern2pt}\renewcommand{\,}{\kern1pt}

\newcommand{\dirac}{/\kern-5pt\partial}

\newcommand{\A}{\mbox{$\widehat A$}}

\def\frac#1#2{{#1\over#2}}
\def\be#1\ee{\begin{equation}#1\end{equation}}

\newcommand{\ot}{\otimes}
\newcommand{\lra}{\longrightarrow}
\renewcommand{\ts}{\textstyle}
0

\newtheorem{theo}{Theorem}[section]
\newtheorem{defi}[theo]{Definition}
\newtheorem{lemma}[theo]{Lemma}
\newtheorem{corol}[theo]{Corollary}
\newtheorem{remark}[theo]{Remark}

\begin{document}
\title{Higher \A-genera on certain non-spin $S^1$-manifolds}

\author{Hayde\'e Herrera\footnote{Department of
Mathematical Sciences, Rutgers University, Camden, NJ 08102,
USA.\hfill\break E-mail: haydeeh@camden.rutgers.edu}
and Rafael
Herrera\footnote{Centro de Investigaci\'on en Matem\'aticas, A. P.
402, Guanajuato, Gto., C.P. 36000, M\'exico. E-mail:
rherrera@cimat.mx} \footnote{Partially supported by CONACYT
grants: J48320-F, J110.387/2006.} }

\date{}

\maketitle

{ \abstract{

We prove the vanishing of higher \A-genera, in the sense of
Browder and Hsiang \cite{Browder}, on smooth $S^1$-manifolds
with finite $\pi_2$ and $\pi_4$.\\


} }

\section{Introduction}

The classical result of Atiyah and Hirzebruch \c{AH} about the
 vanishing of the
\A-genus on Spin manifolds with $S^1$ actions was generalized by
Browder and Hsiang \c{Browder} to higher \A-genera in the
following form.

\begin{theo} {\ns \cite[Theorem 1.8]{Browder}}
Let $M$ be a closed  Spin manifold with a smooth effective action
of a compact, connected, positive-dimensional Lie group $G$. Then
\[p_*([M]\cap \A)=0,\]
where $p\colon M\lra M/G$, and $\A\in H^{4*}(M;\mathbb{Q})$ is the
$\A$ polynomial in the Pontrjagin classes.
\end{theo}

Furthermore, from this theorem they also deduced a higher \A-genus
theorem analogous to Novikov's ``higher signature".    \\

By a closed manifold $M$, we mean a compact manifold without
boundary.   Notice that if $G$ is a compact Lie group not
necessarily connected then we restrict our attention to the
connected component of the identity element.

In this paper, we prove two theorems (Theorems~\ref{theo1}
and~\ref{theo2}) for non-Spin $G$-manifolds
with finite $\pi_2$ and $\pi_4$. They are analogous to those of Browder and
Hsiang \cite{Browder} for Spin manifolds.

\begin{theo}\label{theo1}
Let $M$ be a smooth, closed, connected, oriented
(even-dimensional) $G$-manifold with finite $\pi_2(M)$ and
$\pi_4(M)$, where $G$ is a compact, connected,
positive-dimensional Lie group. Then for any $y\in H^*(M/G,
\mathbb{Q})$ \[(\A \cup p^*(y))[M]=0,\] which implies
\[p_*([M]\cap \A)=0,\]
where $p\colon M\lra M/G$ is the projection map, $\A\in
 H^{4*}(M;\mathbb{Q})$ is the
$\A$ polynomial.
\end{theo}

The proof will make use of the $G$-transversality approach of
Browder and Quinn \cite{Browder-Quinn}, properties of
$G$-transverse submanifolds, and the rigidity of the elliptic
genus on manifolds admitting 2-balanced $S^1$ actions (see below).

\vspace{.2in}

\noi {\it Acknowledgements}.
The authors wish to thank Anand Dessai for fruitful conversations.
The first named author wishes to
thank the Centro de Investigaci\'on en Matem\'aticas and the
Instituto de Matem\'aticas of UNAM (M\'exico) for their
hospitality and support. The second named author wishes to thank
the Max Planck Institute for Mathematics (Bonn) and the
Institut des Hautes \'Etudes Scientifiques
(France) for their hospitality and support.

\section{$S^1$-transverse submanifolds of
manifolds with finite $\pi_2$ and $\pi_4$}

\begin{defi} Let $G$ be a connected Lie group acting smoothly on a
 manifold $M$. Let $H$ be a subgroup of $G$. We denote by
$M^H$ the fixed point set of $H$ on $M$. A $G$-invariant
submanifold $N$ of  $M$ is called {\em transverse} if $N$
intersects $M^H$ transversely for every subgroup $H$ of $G$.
\end{defi}

In order to prove Theorem~\ref{theo1}, we only need to consider a circle action. Thus, we
can choose any circle subgroup $S^1 \subseteq G$.
We denote by $M^{S^1}$
the fixed point set of the circle action. At a fixed point
$p\in M^{S^1}$, the tangent space of $M$
becomes a real representation of $S^1$, whose complexification can be written as
\[T_pM\otimes \mathbb{C}= (t^{m_1}+ t^{-m_1}) + \cdots  + (t^{m_d} + t^{-m_d}) \]
where $t^a$ denotes the representation on which
$\lambda\in S^1$ acts by multiplication by $\lambda^a$, and $d$ is half the
dimension of $M$. The term $(t^{n} + t^{-n})$ corresponds to the representation
\[\lambda=e^{i\theta}\in S^1 \mapsto
\left(\begin{array}{cc}
\cos(n\theta) & -\sin(n\theta) \\
\sin(n\theta) & \cos(n\theta)
\end{array}\right).
\]
The numbers $\pm m_1,\ldots, \pm m_d$ are called the {\rm exponents (or
weights)} of the $S^1$-action at the point $p$.
A circle action is called {\em $2$-balanced} if the parity of
$\sum_{i=1}^dm_i$  does not depend on the connected component of
$M^{S^1}$ (cf. \c{HiBJ}). Since we are only interested in the
parity of $\sum_{i=1}^dm_i$, we do not worry about the choice of
signs.

\begin{lemma}\label{2-balanced} Let $M$ be a $S^1$-manifold
with finite $\pi_2(M)$ and $\pi_4(M)$. Let $N$ be an
$S^1$-transverse submanifold of $M$. Then the $S^1$ action on $N$
is $2$-balanced.
\end{lemma}

{\em Proof}. Since $N$ meets $M^{S^1}$ transversely, for $p\in
N\cap M^{S^1}$
\[T_pM = T_pN + T_pM^{S^1}.\]
Notice that $N^{S^1}= N\cap M^{S^1}$.  Let $p, p'\in N^{S^1}$ lie
in two different components of $N^{S^1}$. The tangent spaces to
$N$ at $p$ and $p'$ become $S^1$ representations so that
\begin{eqnarray}T_pN\otimes \mathbb{C} &=& (t^{n_1(p)} + t^{-n_1(p)}) +\ldots +  (t^{n_k(p)} + t^{-n_k(p)}),\nonumber \\
T_{p'}N\otimes \mathbb{C} &=& (t^{n_1(p')} + t^{-n_1(p')}) +\ldots
+ (t^{n_k(p')} + t^{-n_k(p')}) ,\nonumber\end{eqnarray} and we
have to verify that
\begin{equation}
(n_1(p)+\ldots +n_k(p))-(n_1(p')+\ldots +n_{k'}(p'))\equiv 0 \quad
(\mod 2).\label{difference}
\end{equation}
Observe that the numbers $n_i(p)$ and $n_i(p')$ are, in fact,
exponents of the action of $S^1$ on the manifold
$M$, and that \rf{difference} is the difference of exponents of
the manifold $M$, since the only missing directions of the tangent
space of $M$ are trivial representations (as $N$ and $M^{S^1}$
meet transversely). Since $M$ has finite $\pi_2(M)$ and $\pi_4(M)$,
by \cite[Theorem V]{Bredon}
\[f(t)=T_pM_c - T_{p'}M_c = (1-t)^3 P(t),\]
where $P(t)=\sum b_i t^i$ with only finitely many $b_i$'s
different from  zero. Since real representations are invariant
under the automorphism $t\mapsto t^{-1}$
\[f(t)=f(t^{-1}),\]
i.e.
\[(1-t)^3P(t) = \left( 1-{1\over t}\right)^3 P(t^{-1}).\]
Thus,
\[t^3P(t) = -P(t^{-1}),\]
\[t^3P(t) +P(t^{-1})=0,\]
\[t^3\sum b_i t^i + \sum b_i t^{-i}=0,\]
\[t^{3/2}\sum b_i (t^{i+3/2} + t^{-i-3/2})=0.\]
Since $b_i\not=0$, for every term of the form
$b_i (t^{i+3/2} + t^{-i-3/2})$ there must be another one that cancels it out, i.e.
there must be a $j\not =i$ such that
$b_i=-b_j$ so that either
$j+3/2=i+3/2$ which cannot happen because it contradicts $i\not=j$, or
$-j-3/2=i+3/2$, and $i=-3-j$.
Then, all the terms of $P(t)$ can be grouped according to the corresponding pairs
\[b_it^i + b_jt^j = b_it^i - b_it^{-3-i},\]
which multiplied by $(1-t)^3$ give
\begin{eqnarray}
(b_it^i - b_it^{-3-i})(1-t)^3&=& b_i(t^i+t^{-i}) -3b_i(t^{i+1}+t^{-(i+1)}) \nonumber\\
&+&3b_i(t^{i+2}+t^{-(i+2)}) -b_i(t^{i+3}+t^{-(i+3)}). \nonumber
\end{eqnarray}
Taking the sum of the exponents (with any choice of signs) with multiplicity gives
zero (mod 2),
\[b_i(i) -3b_i(i+1) + 3b_i(i+2) - b_i(i+3) \equiv 0 \quad(\mod 2).\]


\begin{remark}{\rm
Note that the lemma is still valid if we only require the $S^1$
action on $M$ to be $2$-balanced instead of $M$ having
finite $\pi_2(M)$ and $\pi_4(M)$. }\end{remark}

\section{Elliptic genus on manifolds with
$2$-balanced $S^1$-actions}

Let $\ext_c^{\pm}$ be the even and odd complex differential forms
on the oriented, closed, smooth manifold $X$ under the Hodge
$*$-operator, respectively. The signature operator
\[d_s^X =d-*d*\,\,\,\,\colon\,\,\,\, \ext_c^+ \lra \ext_c^-\]
is elliptic and the virtual dimension of its index equals the
signature of $X$, $\sign(X)$. If $W$ is a complex vector bundle on
$X$ endowed with a connection, we can {\em twist} the signature
operator to forms with values in $W$
\[d_s^X\ot W \colon \ext_c^+(W) \lra \ext_c^-(W).\]
This operator is also elliptic and the virtual dimension of its
index is denoted by $\sign(X,W)$.

\begin{defi}
Let $T=TX\otimes\C$ denote the complexified tangent bundle of $X$
and let $R_i$ be the sequence of bundles defined by the formal
series
\[R(q,T)=\sum_{i=0}^{\infty} \,R_i \,\,q^i=
 \bigotimes_{i=1}^{\infty}\ext_{q^i}T
\ot \bigotimes_{j=1}^{\infty}\sym_{q^j}T,\] where
$\sym_tT=\sum_{k=0}^{\infty}\,\sym^kT\,\, t^k$,
$\ext_tT=\sum_{k=0}^{\infty} \,\ext^kT\,\, t^k$, and $\sym^kT$,
$\ext^kT$ denote the $k$-th symmetric and exterior tensor powers
of $T$, respectively. The {\em elliptic genus} of $X$ is defined
as
\begin{equation}
\Phi(X) = \ind(d_s^X\otimes R(q,T))
=\sum_{i=0}^{\infty}
 \sign(X,R_i) \cdot q^i .\label{definition-EllipticGenus}
\end{equation}
\end{defi}

\vspace{.1in}

Note that the first few terms of the sequence $R(q,T)$ are $R_0 =
1$, $R_1 = 2 T$, $R_2 = 2 (T^{\otimes 2} + T)$. In particular, the
constant term of $\Phi(X)$ is $\sign(X)$.

If we assume that $G$ is a group acting on $M$ and commuting with
the elliptic operator, then for $g\in G$ the equivariant index of
$D$ can be defined as \[ \mbox{index}(D)_{G}(g) =
\mbox{trace}(g,\mbox{Ker}D)-\mbox{trace}(g,\mbox{Coker}D).\]

 In an analogous way to the definition of the elliptic genus, now we define the
equivariant elliptic genus with respect to the $S^1$ action by
\begin{equation}
\Phi(X)_{S^1}(\la) =\sum_{i=0}^{\infty}
 \sign(X,R_i)_{S^1}(\lambda)\,\cdot
 q^i,\label{equivariant-elliptic-genus}
\end{equation} where $\lambda\in S^1$.

\begin{theo}\label{main}
Let $X$ be an $2n$-dimensional, oriented, closed, smooth manifold
admitting a smooth $2$-balanced $S^1$-action. Then
\begin{equation}
\Phi(X) = \Phi(X)_{S^1}(\lambda) \label{rig0}
\end{equation}
for every $\lambda\in S^1$.
\end{theo}

{\em Sketch of proof}. The proof of Theorem~\ref{main} is along
the lines of \c{BT}. The equivariant elliptic genus
$\Phi(X)_{S^1}(\lambda)$ turns out to be
 a meromorphic
function on $T_{q^2}=\C^*/q^2$ (the non-zero complex numbers
modulo the multiplicative group generated by $q^2\not =0$). Thus,
the proof of the theorem reduces to proving that
$\Phi(X)_{S^1}(\lambda)$ has no poles at all on $T_{q^2}$, thus
implying that $\Phi(X)_{S^1}(\lambda)$ is constant in $\lambda$.
This follows from applying the Atiyah-Segal equivariant index
theorem and localizing to the $S^1$-fixed point set and other
auxiliary submanifolds. More precisely, one can define the
translate $t_a\Phi(M)_{S^1}(\lambda)$ of $\Phi(M)_{S^1}(\lambda)$
by $a\in \mathbb{C}^*$, to be given by the map at the character
level $\la\mapsto a\la$. In order to prove the rigidity theorem of
$\Phi(M)$, we shall show that none of the translates $t_a\Phi(M)$,
$a\in T_{q^2}$, by points of finite order on $T_{q^2}$, has a pole
on the circle $|\la|=1$. The translates $t_a\Phi(M)$ can be
expressed as {\em twists} of the elliptic genus on some auxiliary
manifolds. The auxiliary submanifolds are the fixed point sets
$X_k$ of the subgroups $\Z_k\subset S^1$, $k\in \mathbb{Z}$. In
doing so, the corresponding expressions have no poles at 1, and
thus $\Phi(X)_{S^1}(\lambda)$ has no poles at points of finite
order in $T_{q^2}$. This argument is valid as long as:

\begin{itemize}
\item[(i)] the submanifolds $X_k$ containing $S^1$-fixed points are
 orientable;
\item[(ii)] it is possible to choose an orientation of $X_k$
compatible with $X$ and all the components $P$ contained in $X_k$.
\end{itemize}

(i) is proved in \cite[Lemma 1]{HH}. (ii) follows as in
\cite[Lemma 9.3]{BT} but using the fact that the action is $2$-balanced.


\begin{corol}\label{A=0}
Let $X$ be a even-dimensional, oriented, closed, connected, smooth
manifold admitting a  $2$-balanced  $S^1$ action. If the $S^1$
action is non-trivial then
\[\A(X) = 0.\]
\end{corol}


The proof follows in the same way as in \cite[Theorem, Section 1.5]{HiSl}. \qd

\begin{corol}\label{vanishing-N}
Let $G$ be a compact positive-dimensional Lie group. Let $N$ be a
compact $G$-transverse submanifold of a connected, oriented,
smooth $G$-manifold with finite $\pi_2(M)$ and $\pi_4(M)$. Then the \A-genus of $N$
vanishes
\[\A(N) =0.\]
\end{corol}

This follows from Lemma~\ref{2-balanced} and Corollary~\ref{A=0}.
\qd

\begin{remark}{\rm
Note that $N$ is not necessarily a Spin manifold since $M$ is not
required to be so. Thus, the vanishing of $\A(N)$ is not a
consequence of the Atiyah-Hirzebruch vanishing theorem. }
\end{remark}

\begin{remark}{\rm
Note that $M$ does not need to be compact for Corollary
\ref{vanishing-N} to hold. }\end{remark}

\section{Vanishing of higher \A-genera}

In this section we provide the proof of Theorem~\ref{theo1}.  The
proof follows that of Theorem 1.8 in \cite{Browder}, and uses the
rigidity of the elliptic genus on manifolds admitting 2-balanced
$S^1$ actions. We also prove Theorem~\ref{theo2}, which can be
thought of as a higher $\hat{A}$-genus theorem for $G$-manifolds
with finite $\pi_2(M)$ and $\pi_4(M)$.\\

{\em Proof of Theorem~\ref{theo1}}.  We use the $G$-transversality
 approach of
Browder and Quinn \cite{Browder-Quinn}.
In it, given a manifold $X$ (not necessarily compact) endowed with
 an action of $G$, they establish a
1-1 correspondence between transverse bordism classes of compact
framed $G$-submanifolds of $X$ of codimension $k$ with
  homotopy classes of maps from $X/G^*$ to
the sphere ${S}^k$, $[X/G^*, S^k]$. Here $X/G^*$ is the 1-point
compactification of $X/G$.

Given a space $Y$, we denote by $\Sigma^t Y$ the $t$-fold reduced
suspension of $Y$, which is homeomorphic to the smash product of
$Y$ and $S^t$, $\Sigma^t Y = Y\wedge S^{t} $.

  To apply \cite[Theorem
4.2]{Browder}, let $y\in H^l(M/G)$. Since rational stable
cohomology and rational stable
 cohomotopy are
isomorphic, we can find $t\in \mathbb{N}$ and a map
\[\rho:\Sigma^t(M/G_{+})\longrightarrow S^{l+t} \]
such that \[\rho^*(g)=\Sigma^t y,\] where $g$ generates
 $H^{l+t}(S^{l+t})$, and $M/G_+$ is the disjoint union of $M/G$
 with a base point.
Notice that $\Sigma^t(M/G_+)=(M/G\times \mathbb{R}^t)^*$,  the one
point compactification of $M/G\times \mathbb{R}^t$. One can
consider $M\times \mathbb{R}^t$ as a $G$-manifold, by extending
the action of $G$ to the $\mathbb{R}^t$ factor by a trivial
action. By \cite[Theorem 4.2]{Browder} there is a compact
transverse framed
 $G$-submanifold $i:N\hookrightarrow M\times \mathbb{R}^t$,
such that
\[p^*\rho^*(g)\cap [(M\times \mathbb{R}^t)^*]
=i_*[N],
\]
i.e. the Poincar\'e dual of $i_*[N]$ is $p^*\rho^*(g)$, which
follows from the construction of the submanifold $N$ in the proof
of Lemma 4.4 in \cite{Browder}.

 Note that
\[\A(N) = i^*\A(M),\]
where $\A(M)\in H^{4*}(M;\mathbb{Q})$. Since $M\times
\mathbb{R}^t$ has finite $\pi_2(M\times
\mathbb{R}^t)$ and $\pi_4(M\times
\mathbb{R}^t)$, and $N$ is a $G$-transverse
 submanifold of $M\times \mathbb{R}^t$,
the $G$ action on $N$ is non-trivial and $\A(N)[N]=0$.  Hence, by
Corollary~\ref{vanishing-N}
\begin{eqnarray}
0&=& \A(N)[N]\nonumber\\ &=& (i^*\A(M))[M]\nonumber\\ &=&
\A(M)(i_*[N])\nonumber\\ &=& \A(M)(p^*\rho^*\cap [M\times
\mathbb{R}^t])\nonumber\\ &=& \A(M)(p^*\Sigma^ty\cap[M\times
\mathbb{R}^t])\nonumber\\ &=& \A(M)(p^*y\cap [M])\nonumber\\ &=&
(\A(M)\cup p^*y)[M].\nonumber
\end{eqnarray}
\qd

Let  $f: M\lra K(\pi_1(M),1)$ be a map, assume that $f_*: \pi_1(M)
\to \pi_1(M)$ is onto, one can  define $\pi'$ to be
$\pi_1(M)/f_*i_*(\pi_1(G))$, where  $i: G\rightarrow M$ is induced
by the action of G on the base point of $M$. Notice that
$i_*(\pi_1(G))$ is contained in the center of $\pi_1(M)$
\cite[page 40]{Hatcher}. Let $\alpha: \pi_1(M) \to \pi'$ be the
projection.

\begin{theo}\label{theo2} Let $M$ be a closed, connected,
smooth manifold with finite $\pi_2(M)$ and $\pi_4(M)$, and let $G$
be
 a compact, connected, positive-dimensional Lie group acting smoothly and effectively on $M$. Let
$f:M\lra K(\pi_1(M),1)$, and $x\in H^*(K(\pi',1);\mathbb{Q})$.
Then $(f^*\alpha^*(x)\cup \A)[M]=0$, where $\A\in
H^*(M;\mathbb{Q})$ is the \A \,\,\,\, polynomial in the Pontrjagin
classes.
\end{theo}

 {\em Proof}. 
By Theorem 1.1 in \c{Browder}, there is a map
$\phi:H_*(M/S^1,\mathbb{Q})\longrightarrow
H_*(K(\pi',1),\mathbb{Q})$
 such that the following diagram
commutes,
\[
\begin{array}{ccc}
H_*(M,\mathbb{Q})&\longrightarrow
\kern-15pt\raise8pt\hbox{\footnotesize{${f_*}$}}
&H_*(K(\pi_1(M),1),\mathbb{Q})\\
\downarrow\hbox{\footnotesize{$p_*$}}&&\downarrow\hbox{\footnotesize{$\alpha_*$}}\\
H_*(M/S^1,\mathbb{Q})&\longrightarrow\kern-15pt\raise8pt\hbox{\footnotesize{${\phi}$}}
& H_*(K(\pi',1),\mathbb{Q}),\\
\end{array}
\] so that $(\A\cup f^*\alpha^*(x))[M] = (\A \cup p^*\phi^*(x) ) [M] = 0$, for every $x\in
H^*(K(\pi',1),\mathbb{Q})$.

\qd

{\small
\renewcommand{\baselinestretch}{0.5}
\newcommand{\bi}{\vspace{-.05in}\bibitem} }

\end{document}